\documentclass[11pt]{article}
\usepackage{amsfonts}
\usepackage{amsfonts}
\usepackage{mathrsfs}
\usepackage{mathrsfs}
\usepackage{mathrsfs}
\usepackage{amsfonts}
\usepackage{amsfonts}
\usepackage{subfigure}
\usepackage{amsthm}
\usepackage[body={16.8cm,24cm}, top=3cm, bottom=3cm]{geometry}

\usepackage{mathrsfs}
\usepackage{amsmath, amssymb}
\usepackage{shortvrb,psfrag}
\usepackage{enumerate}
\usepackage{color}
\usepackage[dvips]{graphicx}

\newtheorem{lemma}{Lemma}[section]

\def\f#1#2{\frac{#1}{#2}}

\newcommand{\be}{\beta}

\newcommand{\ud}{\mathrm{d}}
\newcommand{\ue}{\mathrm{e}}

\newcommand{\ee}{\mathbf{e}}

\numberwithin{equation}{section} 
\def\f#1#2{\frac{#1}{#2}}

\def\<{\langle}
\def\>{\rangle}

\def\be{\begin{equation}}
\def\ee{\end{equation}}
\def\bee{\begin{eqnarray}}
\def\beee{\begin{eqnarray*}}
\def\eee{\end{eqnarray}}
\def\eeee{\end{eqnarray*}}

\begin{document}

\begin{center}
{\bf\Large A note to paper:  
Axial symmetry and classification of stationary solutions of Doi-Onsager 
equation on the sphere with Maier-Saupe potential
}

Hailiang Liu\footnote{Department of Mathematics, Iowa State University, Ames, IA 50011-2064, U.S., hliu@iastate.edu},
\;
 Hui Zhang\footnote{LMCS and School of Mathematical Sciences,
Beijing Normal University, Beijing, 100875, P.R. China,
hzhang@bnu.edu.cn}, 
\; and  \; Pingwen Zhang\footnote{LMAM and
School of Mathematical Sciences, Peking
University, Beijing, 100871, P.R. China, pzhang@pku.edu.cn}
\end{center}
\indent ABSTRACT. This note serves to provide additional details for the proof of Lemma 3.6 
 in our paper 
 [Liu, Zhang and Zhang, Comm. Math. Sci., 3(2005), pp.201-218].    Moreover, we will also present an alternative, yet simpler, proof based on arguments in  [Wang, Zhang and Zhang, CPAM, 68(2015), no. 8, 1326-1398]. 
\\[0.3cm]

In \cite{LZZ05} we showed that in order to determine the number of solutions to the Doi-Onsager equation on the sphere with Maier-Saupe potential,  it suffices to determine the number of zeros of $B(\eta,\alpha)$ in term of the intensity parameter $\alpha>0$, where
\begin{equation}
B(\eta, \alpha)=\frac{3e^{-\eta}}{\int_0^1 e^{-\eta z^2}dz} -\left(
3-2\eta +\frac{4\eta^2}{\alpha} \right).
\end{equation}
Let us state the original Lemma 3.6 in \cite{LZZ05}, and reproduce its proof with supplemental details mainly in steps 2 and 5. 
The equation numbers follow those numbered in \cite{LZZ05}. 
\begin{lemma}\label{lemma-3-4}
The number of zeros of $B(\eta,\alpha)$  is determined by the
intensity  $\alpha$ as follows:

(i). If $\alpha>7.5$, $B(\eta,\alpha)$ has  three zeros
$\eta^*_1<0,\eta^*_2>0$ and $\eta_0=0$.

(ii). If $\alpha=7.5$, $B(\eta,\alpha)$ has  two zeros
$\eta^*_1<0$ and $\eta_0=0$.

(iii). There exists an $\alpha^*\in (20/3,7.5)$ such that
$B(\eta,\alpha)$ has  three zeros $\eta^*_1<0,\eta^*_2<0$ and
$\eta_0=0$ for $\alpha^*<\alpha<7.5$.

(iv). If $ \alpha=\alpha^*$,  $B(\eta,\alpha)$ has  two zeros
$\eta^*_1<0$ and $\eta_0=0$.

 (v). If $0< \alpha< \alpha^*$, $B(\eta,\alpha)$ has  one zero
$\eta_0=0$.
\end{lemma}
 \noindent{\bf Proof. }  This proof is divided into five steps.

 \noindent {\bf Step 1.}  We first show that
 \be
 \label{large}B(\pm M,\alpha)<0 \quad \mbox{ for } M \gg 1.\ee
Note that, for $\eta>0$, the mean value theorem gives  $3e^{-\eta
}/\int_0^{1} e^{-\eta z^2}dz= 3e^{-\eta (1-\gamma^2)}\in (0,3)$
for some $\gamma\in (0,1)$. For $\eta <0$, we have from (3.39)
\beee e^{-\eta}=\int_0^1 e^{-\eta z^2}dz-2\eta \int_0^1
z^2e^{-\eta z^2}dz&\leq  (1-2\eta) \int_0^1 e^{-\eta z^2}dz.\eeee
This  shows  $0<\f{3e^{-\eta }}{\int_0^{1} e^{-\eta z^2}dz}\leq
3(1-2\eta)$. Thus for large $|\eta|\gg 1, B$ is determined by the
quadratic term  $-4\eta ^2/\alpha$. \eqref{large} follows.

{\bf Step 2.} If  $\alpha>7.5$  we will show that $B(\eta,\alpha)$
has at least two zeros $\eta_1<0$, and $\eta_2>0$. We know that
$\eta=0$ is always a zero  since $B(0, \alpha)=0$. Moreover,  let $A_k=\int_{0}^1 z^k e^{-\eta z^2}dz$, then  $\partial_\eta A_k =-A_{k+2}$ for any integer $k\geq 0$. By using (3.40)
we have %
\bee\label{3.22} B_\eta(\eta ,\alpha)&=&3 e^{- \eta
}\f{A_2-A_0}{A_0^2}-\left(  -2 +\f{8\eta}{\alpha}\right).
\eee
This yields
$B_\eta(0,\alpha)=0$. Therefore, $\eta=0$ is a double zero of $B$
for every $\alpha$. Hence the local shape of $B$ 
hinges on the sign of $B_{\eta\eta}(0, \alpha)$. From \eqref{3.22}, we  have
\begin{eqnarray}
\label{3.42}
  B_{\eta\eta}(\eta ,\alpha)&=&\frac{3 e^{-\eta}}{A_0^3}(A_0^2-2A_0A_2-A_0A_4+2A_2^2)-\frac{8}{\alpha}\nonumber\\
  &=&\frac{3 e^{-\eta}}{A_0^3}[A_0(A_0-3A_2)+A_0(A_2-A_4)+2A_2^2]-\frac{8}{\alpha}.
\end{eqnarray}
Therefore a further calculation
from \eqref{3.42} gives
\be
\label{bdd}B_{\eta\eta}(0,\alpha)=\f{16}{15\alpha}\left(\alpha-\f{15}{2}\right).\ee
Thus  $B_{\eta\eta}(0, \alpha)>0$ for $\alpha>7.5$. $B$ is locally
convex near $\eta=0$. This together with \eqref{large} implies
that there exist at least two zeros $\eta^*_1<0$ and $\eta^*_2>0$
besides $\eta=0$.  Now  we assume $\eta^*$ is a zero of $B$, i.e.,
\be\label{zero1}\f{3e^{-\eta^* }}{\int_0^{1} e^{-\eta^*
z^2}dz}=3-2\eta^* +\f{4{\eta^*} ^2}{\alpha} .\ee
In virtue of (3.38), i.e. $A_2=\frac{\alpha-2\eta}{3\alpha}A_0$, and  \eqref{3.22}
we have %
\bee\label{3.45} B_\eta(\eta^* ,\alpha)&=&\f{3e^{-\eta^*  }}{\int_0^{1} e^{-\eta^*
z^2}dz}\left(-1+\f{\alpha-2\eta^*}{3\alpha}\right
)+2\left(1-\f{4\eta^*}{\alpha}\right)\nonumber\\
&=&\f{-2e^{-\eta^*  }}{\int_0^{1} e^{-\eta^*
z^2}dz}\left(1+\f{\eta^*}{\alpha}
\right)+2\left(1-\f{4\eta^*}{\alpha}\right).\eee
Inserting \eqref{zero1} into \eqref{3.45} gives
\bee \label{th-3}B_\eta(\eta^* ,\alpha)&=& -\f{8\eta^*}{3\alpha^2}
\left({\eta^*}^2+\f{\eta^*\alpha}{2}
+\f{\alpha}{2}\left(\f{15}{2}-\alpha\right)\right)\\
 \label{th-4}&=&\left\{\begin{array}{ll} -\f{8\eta^*}{3\alpha^2} (\eta^*-\bar{\eta}_1)
(\eta^*-\bar{\eta}_2),&  \quad \mbox{if}\quad \alpha>20/3,\\
-\f{8\eta^*}{3\alpha^2} (\eta^*+\f{\alpha}{4})^2,&  \quad
\mbox{if}\quad \alpha=20/3,\\
-\f{8\eta^*}{3\alpha^2}
[(\eta^*+\f{\alpha}{4})^2+\f{15\alpha}{4}(1-\f{3\alpha}{20})],&
\quad \mbox{if}\quad \alpha<20/3,\end{array}\right. \eee
where
\be
\label{th-5}\bar{\eta}_1=-\f{\alpha}{4}\left(1+3\sqrt{1-\f{20}{3\alpha}}\right),
\quad
\bar{\eta}_2=-\f{\alpha}{4}\left(1-3\sqrt{1-\f{20}{3\alpha}}\right).\ee
From \eqref{th-4}, we  see that $\eta^*_1\leq \bar{\eta}_1<-\alpha/2$,
and $\eta^*_2\geq\bar{\eta}_2>0$. 

In fact  it is impossible that $\eta^*_1=\bar{\eta}_1$ or $\eta^*_2=\bar{\eta}_2$.
If $\eta^*_1=\bar{\eta}_1$ or $\eta^*_2=\bar{\eta}_2$, then 
$B_\eta(\eta^* ,\alpha) =0$ from \eqref{th-4}; we would have to calculate
$B_{\eta\eta}(\eta^* ,\alpha)$ to reach a contradiction.  For $\eta=\eta^*$, from (3.38) we have $\alpha=\frac{2\eta^* A_0}{A_0-3A_2}$, that gives $A_2/A_0=\frac{\alpha-2\eta^*}{3\alpha}$. Using the fact $B_\eta(\eta^*, \alpha)=0$ 
we also obtain $\alpha=\frac{A_0}{A_2-A_4}$.  Thus
\begin{eqnarray}
B_{\eta\eta}(\eta^*,\alpha) &=&\frac{3 e^{-\eta^*}}{A_0}\left[\frac{A_0-3A_2}{A_0}+\frac{A_2-A_4}{A_0}+2\frac{A_2^2}{A_0^2}\right]-\frac{8}{\alpha}
\nonumber \\
&=&
\frac{3 e^{-\eta^*}}{A_0}\Big(\frac{2\eta^*}{\alpha}+\frac{1}{\alpha}+\frac{2}{9}\frac{(\alpha-2\eta^*)^2}{\alpha^2}\Big)-\frac{8}{\alpha}
\nonumber\\
  \quad &=&\frac1\alpha\Big(3-2\eta^*+\frac{4{\eta^*}^2}{\alpha}\Big)\Big({2\eta^*}+1
+\frac{2}{9}\frac{(\alpha-2\eta^*)^2}{\alpha}\Big)-\frac{8}{\alpha},\nonumber
\end{eqnarray}
where we have used (3.45). Now we denote
$$A={\eta}^2+\f{\eta\alpha}{2}
+\f{\alpha}{2}\left(\f{15}{2}-\alpha\right)=(\eta-\bar{\eta}_1)
(\eta-\bar{\eta}_2).
$$
If $\eta^*=\bar\eta_1$ or $\bar\eta_2$, then $A=0$ at $\eta=\eta^*$. Therefore we get
\begin{eqnarray}
B_{\eta\eta}(\eta^*,\alpha) &=&\frac1\alpha(3-2\eta^*-2\eta^*-15+2\alpha)\nonumber\\
&&\cdot \Big({2\eta^*}+1
+\frac{2\alpha}{9}-\frac{8\eta^*}{9}+\frac{8}{9}\big[-\frac{\eta^*}{2}-\frac12(\frac{15}{2}-\alpha)\big]\Big)-\frac{8}{\alpha}\nonumber\\
&=&\frac2{3\alpha}\big((\alpha-2\eta^*-6)({2\alpha}+2\eta^*-7)-12\big)\nonumber\\
&=&\frac2{3\alpha}\big(2\alpha^2-2\alpha\eta^*-4{\eta^*}^2-19\alpha+2\eta^*+30\big)\nonumber\\ \label{eq:2deriv}
&=&\frac4{3\alpha}\big(\eta^*+15-2\alpha\big).
\end{eqnarray}
Apparently, for $\alpha>7.5$,  $B_{\eta\eta}(\eta_1^* ,\alpha)<0$ as $\eta_1^*<0$.
For $\eta_2^*=\bar\eta_2=-\frac{\alpha}{4}(1-3\sqrt{1-\frac{20}{3\alpha}})$,
it holds (recalling $\alpha>15/2$)
\begin{align*}
  B_{\eta\eta}(\eta_2^* ,\alpha)&= -\frac{1}{3}\Big(1-3\sqrt{1-\frac{20}{3\alpha}}\Big)+\f{8}{3\alpha}\Big(\f{15}{2}-\alpha\Big)=\sqrt{\frac{3\alpha-20}{3\alpha}}\Big(1-\sqrt{\frac{3(3\alpha-20)}{\alpha}}\Big)<0.
\end{align*}
These facts  imply that $B$ is locally concave near $\eta_i^*$ if $\eta_i^*=\bar\eta_i$ for $i=1$ or 2.

On the other hand, if $\eta_2^*=\bar\eta_2$ we know that $B$ is locally convex near $\eta=0$. Note $B_{\eta}(0,\alpha)=B_{\eta}(\eta_2^*,\alpha)=0$. Therefore there is at least one zero point $\eta^*_3\in (0, \bar\eta_2)$ of $B(\eta, \alpha)$
which  satisfies
$$B_\eta(\eta^*_3, \alpha)\leq 0. $$
However, this is impossible by (3.48). Therefore $\eta_2^*>\bar\eta_2$.
Similarly we can get $\eta_1^*<\bar\eta_1$.

{\bf Step 3.} We now  show that $B(\eta,\alpha)$ has at most two
zeros besides $0$ for $\alpha>7.5$. From \eqref{th-4}, $
B_\eta(\eta^* ,\alpha)<0$ for $\eta^*\in (\bar{\eta}_2,\infty)$.
This implies that there is at most one zero of $B$ in
$(0,\infty)$. Otherwise $B_\eta(\eta^* ,\alpha)$ has to be
negative at another zero. Similarly, there exists at most one
$\eta^*\in (-\infty,\bar{\eta}_1)$. The claim in (i) is thus
proved.

 {\bf Step 4.}  We now consider the case  $\alpha=7.5$, for which we  show that there exist  two zeros
$\eta^*_1<0$ and $0$. In this case, $ B_{\eta}(0 ,\alpha)=
B_{\eta\eta}(0 ,\alpha)=0$. In order to see the local shape of $B$
at $\eta=0$, we calculate $B_{\eta\eta\eta}(0 ,\alpha)$.  From (3.43) we have by a careful calculation that 
\begin{eqnarray}
\label{3.51}
  B_{\eta\eta\eta}(\eta ,\alpha)&=&\frac{3 e^{-\eta}}{A_0^4}[6A_2^3-6A_0A_2(A_2+A_4)+A_0^2(3A_2+3A_4+A_6)-A_0^3].
\end{eqnarray}
So 
\be\label{three} B_{\eta\eta\eta}(0 ,\alpha)=-\f{32}{105}<0.\ee
This means
\be\label{eta}\eta B(\eta,\alpha)<0\quad\quad \mbox{for} \quad
|\eta|\ll 1.\ee \eqref{th-3} with $\alpha=7.5$ gives
\be \label{th-6}B_\eta(\eta^* ,\alpha)=-\f{8{\eta^*}^2}{3\alpha^2}
\left(\eta^*+\f{\alpha}{2} \right),\ee
where $\eta^*$ is assumed to be a zero of $B$.  The local behavior
implied from \eqref{eta} and the negative sign of $B_\eta(\eta^*
,\alpha)$ for $\eta^*>0$ shows that no zero of $B$ exists in
$(0,\infty)$. On the other hand, \eqref{eta}, together with
 $B(-M,\alpha)<0$ shows  that there exists at least one zero  in
$(-\infty,0)$. We denote it by $\eta_1^*$. By \eqref{th-6}, we
know  $B$ has no zeros  in
$(-\alpha/2, 0)$ (otherwise consider the one closest to 0). Moreover, $\eta_1^*\neq -\alpha/2$ (otherwise $B_\eta(\eta_1^* ,\alpha)=0$ and by (\ref{eq:2deriv}) $B_{\eta\eta}(\eta_1^* ,\alpha)<0$ which is impossible since  $B$ has no zeros  in
$(-\alpha/2, 0)$).  Thus $\eta_1^*<-\alpha/2$ and $B_\eta(\eta_1^*, \alpha)>0$. Now we
claim $\eta_1^*$ is a unique zero of $B$ in $(-\infty,0)$.
Otherwise, there should appear  at least two more zeros in
$(-\infty,0)$, which is not allowed by \eqref{th-6}. This proves
(ii).

{\bf Step 5.}  We can show that $B(\eta,
\alpha)$ has no zero in $(0,\infty)$ for $\alpha<7.5$. Otherwise, as $B(0, \alpha)=B_\eta(0, \alpha)=0$ and $B_{\eta\eta}(0, \alpha)<0$ (by \eqref{bdd}),  there must be a positive zero of $B(\eta,\alpha)$ satisfying $B_\eta(\eta_1, \alpha)\ge 0$, which is impossible by \eqref{th-3}.
 Moreover, if $\alpha\leq 20/3$, $B(\eta, \alpha)$ even has no zero in
$(-\infty,0)$. This is ensured by the fact $B(-M,\alpha)<0$ and
the sign constrained by \eqref{th-4}, as argued in step 4.

In order to identify the second critical value $\alpha^*\in
(20/3,7.5)$, we need to use a continuity argument.

First for
$7.5-\delta<\alpha<7.5$, $\delta>0$ small, there are at least two zeros
$\eta_1^*,\eta_2^*<0$ of $B$. In fact for this range of $\alpha$,
$B_{\eta\eta}(0,\alpha)<0$. Thus $B(\eta,\alpha)$ is locally concave
near $\eta=0$. We also know that $B(\eta,7.5)>0$ for $\eta\in
(-\delta_1, 0)$ from \eqref{eta}.   This implies that there exists a point $\eta_0\in (-\delta_1, 0)$ such that $B(\eta_0,\alpha)>0$ for $7.5-\delta<\alpha<7.5$ by the continuity of $B$ in $\alpha$.  This together with $B(-M,\alpha)<0$ shows that
there are two zeros of $B$ in $(-\infty,0)$ for $7.5-\delta<\alpha<7.5$.

Secondly, we claim that for $\alpha<7.5$, $B$ has at most two zeros  in $(-\infty,  0)$. This can be concluded from the following facts (see \eqref{th-5} for definitions of $\bar\eta_{1}, \bar\eta_2$, both are negative for $\alpha<7.5$):
  \begin{itemize}
 \item[(1).] $B (\cdot, \alpha)$ has  no zero point in $(\bar\eta_2, 0)$. Otherwise, let $\eta^*\in (\bar\eta_2, 0)$  be the largest zero point. Then $B_{\eta}(\eta^*, \alpha)<0$ which is not allowed by  \eqref{th-4};
\item[(2).] $\bar\eta_2$ is not  a zero point. Otherwise, by \eqref{eq:2deriv} we have $B_{\eta\eta}(\bar\eta_2, \alpha)>0$ which contradicts  to \eqref{th-4};
 \item[(3).]$B (\cdot, \alpha)$ has  at most one zero point in $(\bar\eta_1, \bar\eta_2)$. Otherwise there is a zero point $\eta^*\in (\bar\eta_1, \bar\eta_2)$ satisfying $B_{\eta}(\eta^*, \alpha)\ge 0$ which is impossible by  \eqref{th-4} again;
\item[(4).] $B (\cdot, \alpha)$ has  at most one zero point in $(-\infty, \bar\eta_1)$. Otherwise there is a zero point $\eta^*\in (-\infty, \bar\eta_1)$ satisfying $B_{\eta}(\eta^*, \alpha)\le 0$ which also contradicts to \eqref{th-4};
 \item[(5).] If $\bar\eta_1$ is  a zero point, then $B(\bar\eta_1, \alpha)=B_{\eta}(\bar\eta_1, \alpha)=0, B_{\eta\eta}(\bar\eta_1, \alpha)<0$ (using  \eqref{eq:2deriv}). Thus repeating the argument in (3) or (4), we know there is no zero point of $B$ in $(\bar\eta_1, \bar\eta_2)$ or $(-\infty, \bar\eta_1)$.
  \end{itemize}

 Now, choose $\alpha_0\in (7.5-\delta, 7.5)$ such that $B(\eta, \alpha_0)$ has exactly two zeros $\eta_{0,1}<\eta_{0,2}$ on $(-\infty, 0)$. Then $B(\cdot, \alpha_0)$ is negative on  $(-\infty, \eta_{0,1})\cup ( \eta_{0,2}, 0)$.  As  $B(\eta, \cdot)$ is a monotonically increasing function of $\alpha$, we know that for all $\alpha<\alpha_0$,  $B(\eta, \alpha)<0$ for $\eta \in (-\infty, \eta_{0,1})\cup ( \eta_{0,2}, 0)$. Let
\begin{align}
  \alpha^*=\sup \{\alpha<\alpha_0\,\, | B(\eta,\alpha)<0 \mbox{ for all }\eta \in (\eta_{0,1}, \eta_{0,2})\}.
\end{align}
Then for $\alpha<\alpha^*$, $B$ has no zeros on $(-\infty, 0)$. Now we prove $B$ has two zeros for $\alpha\in(\alpha^*, 7.5)$ and has one zero for $\alpha=\alpha^*$. Thus $\alpha^*\in (20/3, 7.5)$ is
the corresponding critical value.
%

 Apparently, $B(\eta, \alpha^*)$ must have (at least) a zero $\eta^*$ in $(\eta_{0,1}, \eta_{0,2})$. Moreover, $B_\eta(\eta^*, \alpha^*)=0$, and thus,  \eqref{th-4} implies $\eta^*\in\{\bar\eta_1(\alpha^*), \bar\eta_2(\alpha^*)\}$.
However, $\eta^*\neq \bar\eta_2$ by the fact (2), thus $\eta^*=\bar\eta_1(\alpha)$ is the unique zero of $B(\cdot, \alpha^*)$ in $(-\infty, 0)$.
Since $B(\eta^*, \cdot)$ is monotonically increasing function of $\alpha$, we know $B(\eta^*, \alpha)>0$ for all $\alpha>\alpha^*$.
Thus, there are exactly two zeros of $B$ with one of them, say $\eta_1^*$, belongs to $(-\infty, \eta^*)$ and another one $\eta_2^*\in (\eta^*, 0)$ for all $\alpha\in(\alpha^*, 7.5)$.
Using again the fact that $B(\eta^*, \alpha)$ is monotonically increasing function of $\alpha$,
we know $\eta_1^*$ is a decreasing function of $\alpha$, while $\eta_2^*$  is a increasing function of $\alpha$.

These all together finish the proof of (iii)-(v). \hfill{$\Box$}

\vspace{2mm}

{\bf An alternative proof}:

The core part of this proof can be found in Lemma A.1 in \cite{WZZ15}. 

{\bf  Step 1: an equivalent formulation}

Let $A_k(\eta)=\int_{0}^1z^k\ue^{-\eta{z^2}}\ud{z}$. Then the use of integration by parts gives 
$$
 (k+1){A_k}-2\eta A_{k+2}=\ue^{-\eta},\,\,
\text{for }k\ge 0.
$$ 
Consequently, we have
\begin{align}
\nonumber
&\frac{3\ue^{-\eta}}{\int_0^1\ue^{-\eta{z^2}}\ud{z}}=3-2\eta+\frac{4\eta^2}{\alpha}  \quad \Longleftrightarrow\quad
 \frac{3(A_0-2\eta A_2)}{A_0}=3-2\eta +\frac{4\eta^2}{\alpha}\\ 
 \quad \Longleftrightarrow\quad&4\eta^2\Big(\frac{A_2-A_4}{A_0}-\frac{1}{\alpha}\Big)=0
\quad \Longleftrightarrow\quad\eta=0,\quad \text{or} \quad \alpha=\frac{A_0}{A_2-A_4}. 
\end{align}
(This formulation was also shown in the end of \cite{LZZ05}, see the last line of pp. 217)

\medskip

{\bf Step 2: number of solutions}

It suffices to explore the solutions to $\alpha=f(\eta):=\frac{A_0(\eta)}{A_2(\eta)-A_4(\eta)}$. 
Using  $A_k'(\eta)=-A_{k+2}$, we have
\begin{align*}
 f'(\eta)=\Big(\frac{A_0}{A_2-A_4}\Big)'
=\frac{-A_2(A_2-A_4)+A_0(A_4-A_6)}{(A_2-A_4)^2}  .
\end{align*}
Using the fact that
\begin{align*}
&\frac{\partial}{\partial\eta}\Big(\ue^{\eta}\big({A}_0(A_4-A_6)-A_2(A_2-A_4)\big)\Big)\\
&=\frac{1}{2}\frac{\partial}{\partial\eta}\int_{0}^1\int_{0}^1
(x^2y^4+x^4y^2-x^6-y^6+x^4+y^4-2x^2y^2)\ue^{-\eta(x^2+y^2-1)}\ud{x}\ud{y}\\
&=\frac{1}{2}\int_{0}^1\int_{0}^1
(x^2-y^2)^2(1-x^2-y^2)^2\ue^{\eta(x^2+y^2-1)}\ud{x}\ud{y}>0,
\end{align*}
we deduce that ${A}_0(A_4-A_6)-A_2(A_2-A_4)=0$ has only one root $\eta_*$,  hence $f'(\eta)(\eta-\eta^*)>0$ for $\eta\not=\eta^*$. Note that from $f'(0)=5/7>0$ it follows that $\eta^*<0$.  
Hence we have 

{\it Conclusion 1: $f(\eta)$ is monotonically decreasing (increasing) on $(-\infty, \eta_*]$ ($[\eta_*,+\infty)$).}

Thus $f(\eta)$ has a unique global minimizer $\alpha_*=f(\eta_*)>0$. In addition, we have
\begin{align*}
  \frac{\eta}{f(\eta)}=\frac{\eta(A_2-A_4)}{A_0}=\frac{A_0-3A_2}{2A_0}\in (-1, \frac12) \qquad  (\text{since } 0< A_2<A_0 ).
\end{align*}
Therefore $f(\eta)>2\eta$ and $f(\eta)>-\eta$, which implies

{\it Conclusion 2: $f(\eta)\to +\infty$ as $\eta\to\pm\infty$.}

\medskip

Combining the above conclusions, we have that:
\begin{itemize}
  \item  For $\alpha>\alpha_*$, the equation $\alpha=f(\eta)$ has exactly two solutions $\eta_1<\eta_*<\eta_2$;
If $\alpha >(<)7.5$, we have $\eta_2>(<)0$ (since $f(0)=7.5$).
  \item For $\alpha=\alpha_*$, the equation $\alpha=f(\eta)$ has only one solution $\eta=\eta_*$;
  \item For  $\alpha<\alpha_*$, the equation $\alpha=f(\eta)$ has no solution.
\end{itemize}
These yield the conclusions in Lemma 3.6 except for the justification of $\alpha_*>20/3$.

\vspace{2mm}

\noindent {\bf Acknowledgments}: We are very grateful to Professor John Ball for pointing out to us a lack of clarity in some parts of Step 2 and 5 in the original proof of Lemma 3.6. We also thank Wei Wang for his help on the second proof.

\bibliographystyle{abbrv}

\end{document}